# Coverage of space in Boolean models

Rahul Roy[1],*

*Indian Statistical Institute*

**Abstract:** For a marked point process $\{(x_i, S_i)_{i \geq 1}\}$ with $\{x_i \in \Lambda : i \geq 1\}$ being a point process on $\Lambda \subseteq \mathbb{R}^d$ and $\{S_i \subseteq R^d : i \geq 1\}$ being random sets consider the region $C = \cup_{i \geq 1}(x_i + S_i)$. This is the *covered* region obtained from the Boolean model $\{(x_i + S_i) : i \geq 1\}$. The Boolean model is said to be *completely covered* if $\Lambda \subseteq C$ almost surely. If $\Lambda$ is an infinite set such that $\mathbf{s} + \Lambda \subseteq \Lambda$ for all $\mathbf{s} \in \Lambda$ (e.g. the orthant), then the Boolean model is said to be *eventually covered* if $\mathbf{t} + \Lambda \subseteq C$ for some $\mathbf{t}$ almost surely. We discuss the issues of coverage when $\Lambda$ is $\mathbb{R}^d$ and when $\Lambda$ is $[0, \infty)^d$.

## 1. Introduction

A question of interest in geometric probability and stochastic geometry is that of the complete coverage of a given region by smaller random sets. This study was initiated in the late 1950's. An account of the work done during that period may be found in Kendall and Moran (1963). A similar question is that of the connectedness of a random graph when two vertices $u$ and $v$ are connected with a probability $p_{u-v}$ independent of other pairs of vertices. Grimmett, Keane and Marstrand (1984) and Kalikow and Weiss (1988) have shown that barring the 'periodic' cases, the graph is almost surely connected if and only if $\sum_i p_i = \infty$.

Mandelbrot (1972) introduced the terminology *interval processes* to study questions of coverage of the real line $\mathbb{R}$ by random intervals, and Shepp (1972) showed that if $S$ is an inhomogeneous Poisson point process on $\mathbb{R} \times [0, \infty)$ with density measure $\lambda \times \mu$ where $\lambda$ is the Lebesgue measure on the $x$-axis and $\mu$ is a given measure on the $y$-axis, then $\cup_{(x,y) \in S}(x, x+y) = \mathbb{R}$ almost surely if and only if $\int_0^1 dx \exp(\int_x^\infty (y-x)\mu(dy)) = \infty$. Shepp also considered random Cantor sets defined as follows: let $1 \geq t_1 \geq t_2 \geq \ldots$ be a sequence of positive numbers decreasing to 0 and let $\mathcal{P}_1, \mathcal{P}_2, \ldots$ be Poisson point processes on $\mathbb{R}$, each with density $\lambda$. The set $V := \mathbb{R} \setminus (\cup_i \cup_{x \in \mathcal{P}_i} (x, x + t_i))$ is the random Cantor set. He showed that $V$ has Lebesgue measure 0 if and only if $\sum_i t_i = \infty$. Moreover, $P(V = \emptyset) = 0$ or 1 according as $\sum_{n=1}^\infty n^{-2} \exp\{\lambda(t_1 + \cdots + t_n)\}$ converges or diverges.

In recent years the study has been re-initiated in light of its connection to percolation theory. Here we have a marked point process $\{(x_i, S_i)_{i \geq 1}\}$ with $\{x_i : i \geq 1\}$ being a point process on $\Lambda \subseteq \mathbb{R}^d$ and $S_i \subseteq \mathbb{R}^d$ being random sets. Let $C = \cup_{i \geq 1}(x_i + S_i)$ be the *covered* region of the Boolean model $\{(x_i + S_i) : i \geq 1\}$.

The simplest model to consider is the Poisson Boolean model, i.e., the process $\{x_i : i \geq 1\}$ is a stationary Poisson point process of intensity $\lambda$ on $\mathbb{R}^d$ and

$$S_i = [0, \rho_i]^d, \quad i \geq 1 \tag{1}$$

---

*This research is supported in part by a grant from DST.
[1]Statistics and Mathematics Unit Indian Statistical Institute, 7 SJS Sansanwal Marg, New Delhi 110016, India, e-mail: `rahul@isid.ac.in`
*AMS 2000 subject classifications:* primary 05C80, 05C40; secondary 60K35.
*Keywords and phrases:* Poisson process, Boolean model, coverage.

are $d$-dimensional cubes the lengths of whose sides form an independent i.i.d collection $\{\rho_i : i \geq 1\}$ of positive random variables. (Alternately, $S_i$'s are $d$-dimensional spheres with random radius $\rho_i$.)

In this case, Hall (1988) showed that

**Theorem 1.1.** $C = \mathbb{R}^d$ *almost surely if and only if* $E\rho_1^d = \infty$.

More generally, Meester and Roy (1996) obtained

**Theorem 1.2.** *For* $\{x_i : i \geq 1\}$ *a stationary point process, if* $E\rho_1^d = \infty$ *then* $C = \mathbb{R}^d$ *almost surely.*

The above results relate to the question of *complete coverage* of the space $\mathbb{R}^d$.

Another question which arises naturally in the Poisson Boolean model is that of *eventual coverage* (see Athreya, Roy and Sarkar [2004]). Let $\{x_i : i \geq 1\}$ be a stationary Poisson process of intensity $\lambda$ on the orthant $\mathbb{R}_+^d$ and the Boolean model is constructed with random squares $S_i$ as above yielding the covered region $C = \cup_{i \geq 1} [x_i(1), x_i(1) + \rho_i] \times \cdots \times [x_i(d), x_i(d) + \rho_i]$. In that case, $P(\mathbb{R}_+^d \subseteq C) = 0$, however we may say that $\mathbb{R}_+^d$ is eventually covered if there exists $0 < t < \infty$ such that $(t, \infty)^d \subseteq C$. In this case, there is a dichotomy *vis-a-vis* dimensions in the coverage properties. In particular, while eventual coverage depends on the intensity $\lambda$ for $d = 1$, for $d \geq 2$ there is no such dependence.

**Theorem 1.3.** *For* $d = 1$,
(a) *if* $0 < l := \liminf_{x \to \infty} xP(\rho_1 > x) < \infty$ *then there exists* $0 < \lambda_0 \leq \frac{1}{l} < \infty$ *such that*

$$P_\lambda(\mathbb{R}_+ \text{ is eventually covered by } C) = \begin{cases} 0 & \text{if } \lambda < \lambda_0 \\ 1 & \text{if } \lambda > \lambda_0; \end{cases}$$

(b) *if* $0 < L := \limsup_{x \to \infty} xP(\rho_1 > x) < \infty$ *then there exists* $0 < \frac{1}{L} \leq \lambda_1 < \infty$ *such that*

$$P_\lambda(\mathbb{R}_+ \text{ is eventually covered by } C) = \begin{cases} 0 & \text{if } \lambda < \lambda_1 \\ 1 & \text{if } \lambda > \lambda_1; \end{cases}$$

(c) *if* $\lim_{x \to \infty} xP(\rho_1 > x) = \infty$ *then for all* $\lambda > 0$, $\mathbb{R}_+$ *is eventually covered by* $C$ *almost surely* $(P_\lambda)$;
(d) *if* $\lim_{x \to \infty} xP(\rho_1 > x) = 0$ *then for any* $\lambda > 0$, $\mathbb{R}_+$ *is not eventually covered by* $C$ *almost surely* $(P_\lambda)$.

**Theorem 1.4.** *Let* $d \geq 2$, *for all* $\lambda > 0$.
(a) $P_\lambda(\mathbb{R}_+^d \text{ is eventually covered by } C) = 1$ *whenever* $\liminf_{x \to \infty} xP(\rho_1 > x) > 0$.
(b) $P_\lambda(\mathbb{R}_+^d \text{ is eventually covered by } C) = 0$ *whenever* $\lim_{x \to \infty} xP(\rho_1 > x) = 0$.

In 1-dimension for the discrete case we may consider a Markov model as follows: $X_1, X_2, \ldots$ is a $\{0, 1\}$ valued Markov chain and $S_i := [0, \rho_i]$, $i = 1, 2, \ldots$ are i.i.d. intervals where $\rho_i$ is as employed in (1). The region $\cup(i + S_i)1_{X_i=1}$ is the covered region. This model has an interesting application in genomic sequencing (Ewens and Grant [2001]). If $p_{ij} = P(X_{n+1} = j \mid X_n = i)$, $i, j = 0$ or $1$, denote the transition probabilities of the Markov chain then we have

**Theorem 1.5.** *Suppose* $0 < p_{00}, p_{10} < 1$.

(a) If $l := \liminf_{j \to \infty} jP(\rho_1 > j) > 1$, then $P\{C \text{ eventually covers } \mathbb{N}\} = 1$ whenever $\frac{p_{01}}{p_{10}+p_{01}} > 1/l$.

(b) If $L := \limsup_{j \to \infty} jP(\rho_1 > j) < \infty$, then $P\{C \text{ eventually covers } \mathbb{N}\} = 0$ whenever $\frac{p_{01}}{p_{10}+p_{01}} < 1/L$.

Molchanov and Scherbakov (2003) considered the case when the Boolean model is *non-stationary*. For a Poisson point process $\{x_i : i \geq 1\}$ of intensity $\lambda$, we place a $d$-dimensional ball $B(x_i, \rho_i h(\|x_i\|))$ centred at $x_i$ and of radius $\rho_i h(\|x_i\|)$ where $\rho_i$ is as before and $h : [0, \infty) \to (0, \infty)$ is a nondecreasing function. Let

$$C = \cup_{i=1}^{\infty} B(x_i, \rho_i h(\|x_i\|))$$

denote the covered region. Let $\pi_d$ denote the volume of a ball of unit radius in $d$ dimensions and take $h_0(r) = \left(\frac{d}{\lambda \pi_d} \log r\right)^{1/d}$.

**Theorem 1.6.** *Suppose $E\rho^{d+\eta} < \infty$ for some $\eta > 0$ and $h$ is as above.*

(a) *If*

$$l_h := \liminf_{r \to \infty} \left(\frac{h(r)}{h_0(r)}\right)^d > \frac{1}{E(\rho_1^d)} \quad \text{then } P(C = \mathbb{R}^d) > 0,$$

*and*

(b) *if*

$$0 \leq L_h := \limsup_{r \to \infty} \left(\frac{h(r)}{h_0(r)}\right)^d < \frac{1}{E(\rho_1^d)} \quad \text{then } P(C = \mathbb{R}^d) = 0.$$

The result in (a) above cannot be translated into an almost sure result because of the lack of ergodicity in the model.

## 2. Complete coverage

We now sketch the proofs of Theorems 1.1, 1.2 and 1.6.

Let $V = [0,1]^d \setminus C$ denote the 'vacant' region in the unit cube $[0,1]^d$;

$$\begin{aligned} E(\ell(V)) &= E \int_{[0,1]^d} 1_{\{x \text{ is not covered}\}} dx \\ &= \exp(-\lambda E \rho_1^d) \end{aligned} \qquad (2)$$

where $\ell$ stands for the $d$-dimensional Lebesgue measure. Hence, if $E\rho_1^d < \infty$ then $E(\ell(V)) > 0$ and so $P([0,1]^d \subseteq C) < 1$.

Conversely, $E\rho_1^d = \infty$ implies $E(\ell(V)) = 0$ and thus by stationarity, $E(\ell(\mathbb{R}^d \setminus C)) = 0$. Using the convexity of the shapes $S_i$ we may conclude that $P(C = \mathbb{R}^d) = 1$.

Here the Poisson structure was used to obtain the expression (2); for a general process we need to extract, if possible, an ergodic component of the process and show that the Boolean model obtained from this ergodic component covers the entire space when $E\rho_1^d = \infty$. To this end let $\{x_i : i \geq 1\}$ be an ergodic point process with density 1. Let $D_n = [0, 2^{n/d}]^d$ and $E_n = \{$there exists $x_i$ in the annulus $D_{n+1} \setminus D_n$ such that $D_0 \subseteq (x_i + S_i)\}$. Also let $A_m$ be the event that $m$ is the first index such that $\#\{i : x_i \in D_{n+1} \setminus D_n\} \geq a 2^n$ for all $n \geq m$ and for some fixed constant $a$. By ergodicity, $\{A_m : m \geq 1\}$ forms a partition of the probability space and we obtain

$$\begin{aligned}
P(\cap_{k=m}^\infty E_k^c \mid A_m) &\leq P(\cap_{k=m}^\infty \cap_{i:x_i \in D_{k+1}\setminus D_k} \{\rho_i \leq 2^{k/d}+1\} \mid A_m) \\
&\leq \prod_{k=m}^\infty P(\rho_1 \leq 2^{k/d}+1)^{a2^{k+1}} \\
&\leq \prod_{k=m+1}^\infty P(\rho_1^d \leq 2^k)^{a2^k} \\
&\leq \prod_{k=m+1}^\infty \left(\prod_{j=1}^{2^k-1} P(\rho_1^d \leq q2^k + j)\right)^a \\
&= \left[\prod_{k=2^m}^\infty \left(1 - P(\rho_1^d > k)\right)\right]^a \\
&= 0 \text{ if and only if } \sum_{k=2^m}^\infty P(\rho_1^d > k) = \infty.
\end{aligned}$$

This completes the proof of Theorems 1.1 and 1.2.

To prove Theorem 1.6 (a) we study the case when $h(r) = l_h^{1/d} h_0(r)$ and $l_h E\rho_1^d > 1 + \delta$ for some $\delta > 0$. It may be easily seen that for $z \in \mathbb{Z}^d$, $P\{z + (-1/2, 1/2]^d \not\subseteq C\} \leq \exp\{-\lambda\mu(R_z)\}$, where $R_z = \{(r,x) \in [0,\infty) \times \mathbb{R}^d : z + (-1/2, 1/2]^d \subseteq B(x, rh(x))\}$ and $\mu$ is the product measure of the measure governing $\rho_1$ and Lebesgue measure. From the properties of $h_0$ it may be seen, after some calculations, that given $\epsilon > 0$, there exists $r$ such that for $||z|| > r$, $\mu(R_z) \geq (1-\epsilon)\pi_d(h(z))^d E\rho_1^d$. Thus we obtain, for some constant $K$,

$$\begin{aligned}
\sum_{z \in \mathbb{Z}^d} P\{z + (-1/2, 1/2]^d \not\subseteq C\} &\leq K \sum_{z \in \mathbb{Z}^d} \exp\{-d(1+\delta)(1-\epsilon)\log(||z||)\} \\
&= K \sum_{z \in \mathbb{Z}^d} ||z||^{-d(1+\delta)(1-\epsilon)} < \infty
\end{aligned}$$

whenever $(1+\delta)(1-\epsilon) > 1$. Invoking the Borel-Cantelli lemma we have that $z + (-1/2, 1/2]^d \not\subseteq C$ occurs for only finitely many $z \in \mathbb{Z}^d$. Using this we now complete the proof of Theorem 1.6 (a).

The proof of Theorem 1.6(b) is more delicate and we just present the idea here. For $d \geq 2$, if we place points in a spherical shell of radius $n^\gamma$, such that the interpoint distances are maximum and are of the order of $n^\beta$ where $0 < \beta < \gamma$ and $\gamma > 1$ then the number of points one can place on this shell is of the order of $n^{(\gamma-\beta)(d-1)}$. Let $V_n$ be the event that one such point is not covered by $C$. It may be shown that there is a choice of $\gamma$ and $\beta$ such that infinitely many events $V_n$ occur with probability 1. For $d = 1$, the same idea may be used and, in fact, the proof is much simpler.

## 3. Eventual coverage

We discretise the space $\mathbb{R}_+^d$ by partitioning it into unit cells $\{(i_1, \ldots, i_d) + (0, 1]^d : i_1, \ldots, i_d = 0, 1, \ldots\}$ and call a vertex $\mathbf{i} := (i_1, \ldots, i_d)$ green if $x \in \mathbf{i} + (0, 1]^d$ for some point $x$ of the point process. Consider two independent i.i.d. collections of random variables $\{\rho_\mathbf{i}^u\}$ and $\{\rho_\mathbf{i}^l\}$ where the distribution of $\rho_\mathbf{i}^u$ and $\rho_\mathbf{i}^l$ are identical to that of $2 + \lfloor \max\{\rho_1, \ldots, \rho_N\} \rfloor$ and $\max\{0, \lfloor \max\{\rho_1, \ldots, \rho_N\} \rfloor - 1\}$ respectively; here $N$ is an independent Poisson random variable with mean $\lambda$ conditioned to be 1 or more. We now define two discrete models, an *upper* model and a *lower* model: in both these models a vertex $\mathbf{i}$ is open or closed independently of other vertices, and the covered region for the upper model is $\cup_{\{\mathbf{i} \text{ open}\}}(\mathbf{i} + [0, \rho_\mathbf{i}^u])$, and that for the lower model is $\cup_{\{\mathbf{i} \text{ open}\}}(\mathbf{i} + [0, \rho_\mathbf{i}^l])$. Observe that the eventual coverage of the Poisson model ensures the same for the upper model and eventual coverage of the lower model ensures the same for the Poisson model. Thus it suffices to consider the eventual coverage question for a discrete model, as in Proposition 3.1 below, and check that the random variables $\rho^u$ and $\rho^l$ satisfy the conditions of the proposition.

We take $\{X_\mathbf{i} : \mathbf{i} \in \mathbb{N}^d\}$ to be an i.i.d. collection of $\{0,1\}$ valued random variables with $p := P(X_\mathbf{i} = 1)$ and $\{\rho_\mathbf{i} : \mathbf{i} \in \mathbb{N}^d\}$ to be another i.i.d. collection of positive integer valued random variables with distribution function $F(= 1 - G)$ and independent of $\{X_\mathbf{i} : \mathbf{i} \in \mathbb{N}^d\}$. Let $C := \cup_{\{\mathbf{i} \in \mathbb{N}^d \mid X_\mathbf{i} = 1\}}(\mathbf{i} + [0, \rho_\mathbf{i}]^d)$. We first consider eventual coverage of $\mathbb{N}^d$ by $X_\mathbf{i}$.

**Proposition 3.1.** *Let $d \geq 2$ and $0 < p < 1$.*

(a) *if $\lim_{j \to \infty} jG(j) = 0$ then $P_p(C \text{ eventually covers } \mathbb{N}^d) = 0$,*
(b) *if $\liminf_{j \to \infty} jG(j) > 0$ then $P_p(C \text{ eventually covers } \mathbb{N}^d) = 1$.*

We sketch the proof for $d = 2$. For $i, j \in \mathbb{N}$ let $A(i,j) := \{(i,j) \notin C\}$. Clearly,

$$P(A(k,j) \cap A(i,j)) = P(A(k-i,j))P(A(i,j)) \text{ for } k \geq i,$$

i.e., for each fixed $j$ the event $A(i,j)$ is a renewal event. Thus, if, for every $j \geq 1$, $\sum_{i=1}^{\infty} P(A(i,j)) = \infty$ then, on every line $\{y = j\}$, $j \geq 1$, we have infinitely many $i$'s for which $(i,j)$ is uncovered with probability one and hence $\mathbb{N}^d$ can never be eventually covered.

To calculate $P_p(A(i,j))$ we divide the rectangle $[1,i] \times [1,j]$ as in Figure 1. For any point $(k,l)$, $1 \leq k \leq i - j$ and $1 \leq l \leq j$, in the shaded region of Figure 1, we ensure that either $X_{(k,l)} = 0$ or $\rho_{(k,l)} \leq k + j - 1$. The remaining square region in Figure 1 is decomposed into $j$ sub squares of length $t$, $1 \leq t \leq j - 1$ and we ensure that for each point $(k,l)$ on the section of the boundary of the sub square $t$ given by the dotted lines either $X_{(k,l)} = 0$ or $\rho_{(k,l)} \leq t$. So,

$$\begin{aligned} P_p(A(i,j)) &= (1-p)\prod_{t=1}^{j-1}(1 - p + pF(t-1))^{2t+1} \prod_{k=1}^{i-j}(1 - p + pF(k+j-1))^j \\ &= (1-p)\prod_{t=1}^{j-1}(1 - pG(t))^{2t+1} \prod_{k=j+1}^{i}(1 - pG(k))^j. \end{aligned} \quad (3)$$

Now choose $\epsilon > 0$ such that $pj\epsilon < 1$ and get $N$ such that, for all $i \geq N$, $iG(i) < \epsilon$. Taking $c_j := \prod_{t=1}^{j-1}(1 - pG(t))^{2t+1}$ from (3) we have that

$$\begin{aligned} \sum_{i=N}^{\infty} P_p(A(i,j)) &= (1-p)c_j \sum_{i=N}^{\infty} \prod_{k=1}^{i-j}(1 - pG(k+j))^j \\ &= (1-p)c_j \sum_{i=N}^{\infty} e_i \quad \text{(say)}. \end{aligned} \quad (4)$$

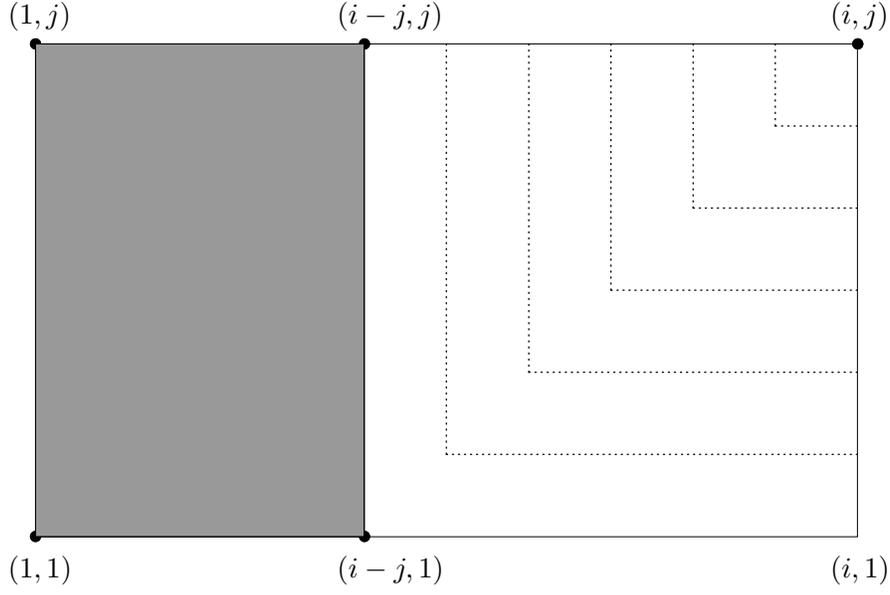

Fig 1. *Division of the rectangle formed by $[1,i] \times [1,j]$.*

For $m \geq N$ we have

$$\begin{aligned}
\frac{e_{m+1}}{e_m} &= (1 - pG(m+1))^j \\
&\geq \left(1 - \frac{p\epsilon}{m+1}\right)^j \\
&= 1 - j\frac{p\epsilon}{m+1} + \sum_{k=2}^{j}(-p)^k \binom{j}{k}\frac{\epsilon^k}{(m+1)^k} \\
&= 1 - \frac{pj\epsilon}{m+1} + \frac{g(m,p,j,\epsilon)}{(m+1)^2},
\end{aligned}$$

for some function $g(m,p,j,\epsilon)$ bounded in $m$. Thus by Gauss' test, as $pj\epsilon < 1$ we have $\sum_{i=N}^{\infty} e_i = \infty$ and hence $\sum_{i=1}^{\infty} P_p(A(i,j)) = \infty$. This completes the proof of the first part of the proposition.

For the next part we fix $\eta > 0$ such that $\eta < \liminf_{j\to\infty} jG(j)$ and get $N_1$ such that for all $i \geq N_1$ we have $iG(i) > \eta$. Also, fix $0 < p < 1$ and choose $a$ such that $0 < \exp(-p\eta) < a < 1$. Let $N_2$ be such that for all $j \geq N_2$ we have $(1 - p\eta j^{-1})^j < a$. For $N := \max\{N_1, N_2\}$, let $i, j \in \mathbb{N}$ be such that $j \geq N$ and $i > j$. Define $A(i,j) := \{(i,j) \notin C\}$. As in (3) we have

$$P_p(A(i,j)) = (1-p) \prod_{k=1}^{i-j}(1 - pG(j+k))^j \prod_{t=1}^{j-1}(1 - pG(t))^{2t+1}. \tag{5}$$

Taking $c_j := \prod_{t=1}^{j-1}(1-pG(t))^{2t+1}$, we have from (5) and our choice of $j$,

$$\sum_{i=N}^{\infty} P_p(A(i,j)) = (1-p)c_j \sum_{i=N}^{\infty} \prod_{k=1}^{i-j}(1-pG(k+j))^j$$

$$= (1-p)c_j \sum_{i=N}^{\infty} b_i \text{ (say)}.$$

For $m \geq N$

$$\frac{b_{m+1}}{b_m} = (1-pG(m+1))^j$$

$$\leq \left(1-p\frac{\eta}{m+1}\right)^j$$

$$= 1 - \frac{pj\eta}{m+1} + \frac{h(m,p,j,\eta)}{(m+1)^2} \quad (6)$$

for some function $h(m,p,j,\eta)$ bounded in $m$; thus by Gauss' test, if $pj\eta > 1$ then $\sum_{i=N}^{\infty} b_i < \infty$ and hence $\sum_{i=1}^{\infty} P_p(A(i,j)) < \infty$.

Now, for a given $p$, let $j' := \sup\{j : pj\eta < 1\}$ and $j_0 := \max\{j'+1, N\}$. We next show that the region $Q_{j_0} := \{(i_1, i_2) \in \mathbb{N}^d : i_1, i_2 \geq j_0\}$ has at most finitely many points that are not covered by $C$ almost surely; there by proving that $C$ eventually covers $\mathbb{N}^d$. For this we apply Borel-Cantelli lemma after showing that $\sum_{(i_1,i_2) \in Q_{j_0}} P_p(A(i_1, i_2)) < \infty$.

Towards this end we have

$$\sum_{i_1, i_2 \geq j_0} P_p(A(i_1, i_2))$$

$$= 2(1-p)\sum_{m=1}^{\infty}\left(\prod_{t=1}^{j_0+m-1}(1-pG(t))^{2t+1}\right.$$

$$\left. \times \sum_{k=m+1}^{\infty}\prod_{i=1}^{k-m}(1-pG(j_0+m+i))^{j_0+m}\right)$$

$$+\sum_{k=1}^{\infty}\prod_{t=1}^{j_0+k-1}(1-pG(t))^{2t+1}.$$

Observe that

$$\sigma_m := \sum_{k=m+1}^{\infty}\prod_{i=1}^{k-m}(1-pG(j_0+m+i))^{j_0+m}$$

$$\leq \sum_{s=1}^{\infty}\prod_{i=1}^{s}\left(1-\frac{p\eta}{j_0+m+s}\right)^{j_0+m},$$

hence as in (6) and the subsequent application of Gauss' test, we have that, for every $m \geq 1$, $\sigma_m < \infty$.

Now let $\gamma_m := \prod_{t=1}^{j_0+m-1}(1-pG(t))^{2t+1}\sigma_m$. Note that an application of the ratio test yields $\sum_{m=1}^{\infty} \gamma_m < \infty$; indeed from (7),

$$\frac{\gamma_{m+1}}{\gamma_m}$$

$$= (1 - pG(j_0 + m))^{2j_0 + 2m + 1} \frac{\sum_{s=1}^{\infty} \prod_{i=1}^{s} (1 - pG(j_0 + m + 1 + i))^{j_0 + m + 1}}{\sum_{s=1}^{\infty} \prod_{i=1}^{s} (1 - pG(j_0 + m + i))^{j_0 + m}}$$

$$= \frac{(1 - pG(j_0 + m))^{2j_0 + 2m + 1}}{(1 - pG(j_0 + m + 1))^{j_0 + m}} \frac{\sum_{s=1}^{\infty} \prod_{i=1}^{s} (1 - pG(j_0 + m + 1 + i))^{j_0 + m + 1}}{1 + \sum_{s=2}^{\infty} \prod_{i=2}^{s} (1 - pG(j_0 + m + i))^{j_0 + m}}$$

$$\leq (1 - pG(j_0 + m))^{j_0 + m + 1} \frac{\sum_{s=1}^{\infty} \prod_{i=1}^{s} (1 - pG(j_0 + m + 1 + i))^{j_0 + m + 1}}{1 + \sum_{s=1}^{\infty} \prod_{i=1}^{s} (1 - pG(j_0 + m + 1 + i))^{j_0 + m}}.$$

Since $\sigma_m < \infty$ for all $m \geq 1$, both the numerator and the denominator in the fraction above are finite. Moreover, each term in the sum of the numerator is less than the corresponding term in the sum of the denominator; yielding that the fraction is at most 1. Hence, for $0 < a < 1$ as chosen earlier

$$\frac{\gamma_{m+1}}{\gamma_m} \leq (1 - pG(j_0 + m))^{j_0 + m + 1}$$

$$\leq a.$$

This shows that $\sum_{m=1}^{\infty} \gamma_m < \infty$ and completes the proof of part (b) of the proposition.

It may now be seen easily that $\rho^u$ and $\rho^l$ satisfy the conditions of Proposition 3.1 and thus Theorem 1.4 holds.

## 4. Markov Model

The relation between the Poisson model and the discrete model explained in Section 3 shows that Theorem 1.3 would follow once we establish Theorem 1.5. In the setup of the Theorem 1.5, for each $k \in \mathbb{N}$ let $A_k := \{k \notin C\}$. To prove Theorem 1.5(a) we show that $\sum_k P(A_k) < \infty$ and an application of the Borel-Cantelli lemma yields the result, while to prove Theorem 1.5 (b) we show that $\sum_k P(A_k) = \infty$. However, the $A_k$'s are not independent and hence Borel-Cantelli lemma cannot be applied. Nonetheless using the Markov property one can show that $P(A_k \cap A_i) = P(A_{k-i})P(A_i)$ and therefore, $A_i$'s are renewal events; so by the renewal theorem, if $\sum_{i=1}^{\infty} P(A_i) = \infty$ then $A_i$ occurs for infinitely many $i$'s with probability one.

For $k \geq 1$, let $P_0(A_k) = P(A_k \mid X_1 = 0)$ and $P_1(A_k) = P(A_k \mid X_1 = 1)$. The following recurrence relations may be easily verified

$$P_0(A_{k+1}) = p_{00} P_0(A_k) + p_{01} P_1(A_k) \tag{7}$$
$$P_1(A_{k+1}) = F(k-1)\left[p_{10} P_0(A_k) + p_{11} P_1(A_k)\right]. \tag{8}$$

We use this to prove Theorem 1.5(b) first. Let $\Psi_0(s) = \sum_{k=k_0}^{\infty} P_0(A_k) s^k$ and $\Psi_1(s) = \sum_{k=k_0}^{\infty} P_1(A_k) s^k$ denote the generating functions of the sequences $\{P_0(A_k) : k \geq k_0\}$ and $\{P_1(A_k) : k \geq k_0\}$ respectively, where $k_0$ is such that for a given $\epsilon > 0$ and $C = L + \epsilon > 0$ (where $L$ is as in the statement of the theorem), $k_0 + (1-C) > 0$, $P_0(A_{k_0}) > 0$, $P_1(A_{k_0}) > 0$, and $F(k-1) \geq 1 - \frac{C}{k+1}$ for $k \geq k_0$. Such a $k_0$ exists by the condition of the theorem.

Using the recurrence relations (7) and (8) we obtain

$$\Psi_1'(s) P(s) \geq Q(s) B(s) + R(s), \tag{9}$$

where

$$\begin{aligned} P(s) &= (1-p_{00}s)^2(1-p_{11}s) + p_{10}s(1-p_{00}s)p_{01}s \\ &= (1-p_{00}s)(1-s)(1-s(1-p_{01}-p_{10})) \\ Q(s) &= (1-p_{00}s)^2(1-C)p_{11} + (1-C)p_{10}p_{01}s(1-p_{00}s) \\ &\quad + p_{10}sp_{01}(1-p_{00}s) + p_{10}p_{00}p_{01}s^2 \\ R(s) &= (1-p_{00}s)^2 k_0 s^{k_0-1} P_1(A_{k_0}) + (k_0+1-C)p_{10}s^{k_0}(1-p_{00}s)P_0(A_{k_0}) \\ &\quad + p_{10}s^{k_0+1}p_{00}P_0(A_{k_0}). \end{aligned}$$

From (9) we have for any $0 < t < 1$

$$\Psi_1(t) \geq e^{\int_0^t \frac{Q(s)}{P(s)}ds} \int_0^t e^{\int_0^s \frac{-Q(r)}{P(r)}dr} \frac{R(s)}{P(s)}ds.$$

Now for $s < 1$, $\frac{Q(s)}{P(s)} = \frac{D}{1-p_{00}s} + \frac{E}{1-s} + \frac{F}{1-s(1-p_{01}-p_{10})}$, for some real numbers $D, E, F$. It may now be seen that $\Psi_1(1) = \infty$ whenever $E > 0$. Also the recurrence relations show that $\Psi_0(1) = \infty$ whenever $\Psi_1(1) = \infty$. A simple calculation now yields that $E > 0$ if and only if $\frac{p_{01}}{p_{10}+p_{01}} < \frac{1}{C} = \frac{1}{L+\epsilon}$. Since $\epsilon$ is arbitrary, we obtain Theorem 1.5(b).

The proof of Theorem 1.5(a) is similar.